\documentclass{article}
\usepackage[T2A]{fontenc}
\usepackage[cp1251]{inputenc}
\usepackage[english,russian]{babel}
\usepackage[tbtags]{amsmath}
\usepackage{amsfonts,amssymb,mathrsfs,amscd}
\usepackage[hyper]{a}
\JournalName{}
\numberwithin{equation}{section}
\theoremstyle{plain}
\newtheorem{theorem}{Theorem}[section]
\newtheorem{lemma}{Lemma}[section]

\theoremstyle{definition}
\newtheorem{definition}{Definition}[section]

\newtheorem{remark}{Remark}[section]

\makeatletter
\renewenvironment{thebibliography}[1]
{\section*{Bibliography}
  \list{\@biblabel{\@arabic\c@enumiv}}%
  {\settowidth\labelwidth{\@biblabel{#1}}%
    \leftmargin\labelwidth
    \advance\leftmargin\labelsep
    \setlength{\itemsep}{0pt}
    \@openbib@code
    \usecounter{enumiv}%
    \let\p@enumiv\@empty
    \renewcommand\theenumiv{\@arabic\c@enumiv}}%
  \sloppy
  \clubpenalty4000
  \@clubpenalty \clubpenalty
  \widowpenalty4000%
  \sfcode`\.\@m}
\makeatother

\begin{document}

\centerline {\bf On the inviscid limit of stationary measures for the stochastic system}
\centerline {\bf of the  Lorenz model for a baroclinic atmosphere}
\vskip10pt
\author
{Yu.\,Yu.~Klevtsova}

\maketitle
\markright{On the inviscid limit of stationary measures}

\footnotetext[0]{This research was carried out within the framework of the Research Plan of Federal
Service for Hydrometeorology and Environmental Monitoring for 2021–2024, the topic 1.1.3
“Development and improvement of the new generation of COSMO-Ru ultra-high-resolution
short-range weather prediction system (with grid steps up to 1 km) based on the ICON
seamless non-hydrostatic atmospheric model”.}

\begin{fulltext}

The paper is concerned with a~nonlinear system of partial differential
equations with parameters and the random external force. This system describes the two-layer
quasi-solenoidal Lorenz model for a~baroclinic atmosphere on a~rotating
two-dimen\-sional sphere.
The stationary measures for the Markov semigroup defined by the solutions of the Cauchy problem for this problem is considered.  One parameter of the system is highlighted -- the coefficient of kinematic viscosity.  The sufficient conditions on the random right-hand side and the other parameters are derived for the existence of a limiting nontrivial point for any sequence of the stationary measures for this system when any sequence of the kinematic viscosity coefficients goes to zero. As it is well known, this coefficient in practice is extremely small.
A number of integral properties are proved for the limiting measure. In addition, these results are obtained for one similar baroclinic atmosphere system.

Bibliography: 21 titles.

{\bf Keywords:}
baroclinic atmosphere, Lorenz model, random external force,  statio\-nary measure, inviscid limit.

\section{Introduction}
\label{s1}

We shall be concerned with the two-layer
quasi-solenoidal model for a~baroclinic atmosphere  in isobaric coordinates,
in which the pressure~$p$ is used in place of the independent vertical coordinate (briefly, the $p$-coordinate system).
The atmosphere is split into two layers as follows: in the first layer the pressure ranges
between $0$ and $500$~mb, and in the second layer, between $500$ and $1000$~mb.
There are five pressure levels: the zeroth level ($p=0$~mb), the first level ($p=250$~mb),
the second layer ($p=500$~mb), the third layer ($p=750$~mb), and the fourth layer ($p=1000$~mb).
An unknown vector function  $u=(u_1,\allowbreak u_2)^{\mathrm T}$ satisfies the following system of equations:
\begin{equation}
\label{eq1.1}
\frac{\partial}{\partial t}A_1 u+\nu A_2u+A_3u+B(u)=\eta,
\qquad t>0.
\end{equation}
The variables $u_1$ and $u_2$  are
the half-sum (the barotropic component) and the half-difference (the baroclinic component), respectively,
of the dimensionless stream function for the divergence-free component
of the velocity vector of horizontal motion in the $p$-coordinate system measured in the first and third levels.
The system~\eqref{eq1.1} is defined on the rotating two-dimensional unit sphere $S$ centred at the origin of the spherical polar
system of coordinates
$$
(\lambda,\varphi),
\qquad
  \lambda \in [0, 2\pi),
\quad
 \varphi \in\biggl[-\frac{\pi}{2},\frac{\pi}{2}\biggr],
 \quad \mu =\sin \varphi , \quad x=(\lambda,\mu).
$$
In~\eqref{eq1.1} vector function $u=u^\omega (t,x)$,
$$
\eta ^\omega (t, x)=\bigl(\eta_1^\omega (t,x),\eta_2^\omega (t,x)\bigr)^{\mathrm T}
$$
is white noise with respect to~$t$,
$\omega \in \Omega $, $(\Omega, F, {\mathbb P})$~ is
a~complete probability space.
In the system~\eqref{eq1.1} numerical parameter $\nu>0$  is the kinematic viscosity coefficient.
The operators in equation~\eqref{eq1.1}
are defined on the sphere~$S$ and are given by
\begin{gather*}
A_1=\begin{pmatrix}
-\Delta & 0
\\
0 & -\Delta+\gamma 	I
\end{pmatrix},
\qquad
A_2=\begin{pmatrix}
\Delta ^2 & 0
\\
0 & \Delta ^2	
\end{pmatrix},
\\
A_3=\begin{pmatrix}
-k_0\Delta & 2k_0\Delta
\\
k_0\Delta & -(2k_0+k_1+\nu \gamma)\Delta+\rho 	I
\end{pmatrix},
\\
B(u)=\bigl(J(\Delta u_1+2\mu,u_1)+J (\Delta u_2, u_2),\,
J(\Delta u_2 -\gamma u_2, u_1)+J(\Delta u_1+2\mu, u_2)\bigr)^{\mathrm T}.
\end{gather*}
Here $B$ is a nonlinear operator, $\gamma,\rho, k_0, k_1\ge 0$ are numerical parameters, $I$ is the identity operator,
$$
J(\psi,\theta)=\psi_{\lambda}\theta_{\mu} - \psi_{\mu}\theta_{\lambda}
$$
is the Jacobian and
$$
\Delta \psi=((1-\mu ^2)\psi_\mu)_\mu+(1 - \mu ^2)^{-1} \psi_{\lambda \lambda}
$$
is the Laplace-Beltrami operator.
The numerical parameters $k_0$ and $k_1$ are proportio\-nal
to the skin friction coefficient and the coefficient of internal interlayer friction, respectively;
the parameter~$\gamma$ is inversely proportio\-nal to the statistical stability parameter of the layer $250$--$750$~mb and $\rho=2 \gamma h$, where $h$~is the parameter responsible for the atmospheric heating from the Earth's surface.

Two is the least number of layers considered by meteorologists in models for a~baroclinic atmosphere.
The number of layers can be arbitrary;
 the larger the number of layers, the more accurate the model.
The barotropy assumes that air density depends only on the pressure. In the baroclinic case, it also depends on other parameters, in our case, on the temperature. The baroclinicity is one of the reasons for the development of vorticity in the atmosphere. The real atmosphere is baroclinic.
E. N. Lorenz~\cite{0} was the first who proposed to study a two-layer quasi-solenoidal model of a baroclinic atmosphere~\eqref{eq1.1} in order to study the regimes of the general circulation of the atmosphere.
This model describes the large-scale dynamics of the atmosphere.
Unlike more complex models, in  case of white noise perturbation    it allows one to obtain the rigorous results  on the existence of a unique stationary measure and on the estimate for the rate of convergence of the distributions of all solutions   of the system~\eqref{eq1.1} lying in a certain class to this measure as $t\to +\infty $ (see~\cite{1}--\cite{52}).
We note that the stationary measure can be understood as a statistical equilibrium, to which the characteristics of all solutions converge as $t\to +\infty $.
In this paper we continue the investigation we began in~\cite{1}--\cite{52}
and examine the inviscid limit of stationary measures.
 As it is well known, in practice, the kinematic viscosity coefficient is extremely small, so meteorologists are willing to replace it with zero. Is such a step justified? Until now, no such mathematical studies have been carried out for baroclinic models. This work is the first step in this direction.

Let
\begin{equation}
\label{eq21.30}
\eta=\nu^\alpha \eta', \qquad k_0=\nu k'_0, \qquad k_1 = \nu k'_1, \qquad \rho = \nu \rho',
\end{equation}
where
$\eta'$, $k'_0$, $k'_1$, $\rho'$ are independent of $\nu$, $\alpha$ is an arbitrary real number. Let  $\delta _0$ denote the Dirac measure. The main result of this paper is the following Theorem.
\begin{theorem}
Let  $m_\nu$ denote the stationary measure for the system~\eqref{eq1.1}  for  $\nu >0$.
Under some condition on
$\eta '\ne 0$  and  for any  $\gamma,\rho', k'_0, k'_1\ge 0$ satisfying the inequality
\begin{equation}
\label{eq7.10}
k'_0\le \min\left\{4k'_1, \frac{4(2+\gamma)}{(2-\gamma)^2}(2k'_1+\rho')\right\},
\end{equation}
the following statements hold.

(i) Let $\alpha = 0.5$.
Then there exists a limiting point $m \ne \delta _0 $ for any sequence of stationary measures
$\{m_{\nu_n}\}_{n=1}^{\infty}$,
$\nu_n\underset{n \to +\infty}{\longrightarrow} 0$,
in the sense of weak-$*$ convergence in some family of probability measures ${\mathscr P}$.
The limiting measure $m$ has some integral properties.

(ii) Let $\alpha > 0.5$.
Then the probability measure  $m_{\nu}$  converges to the Dirac measure $\delta_0$  in the sense of weak-$*$ convergence in ${\mathscr P}$ as $\nu\to 0$.

(iii) Let $\alpha < 0.5$.
Then the set of probability measures   $\{m_{\nu}\}_{\nu >0}$ has no accumula\-tion points  in the sense of weak-$*$ convergence in ${\mathscr P}$ as $\nu\to 0$.

\end{theorem}

Thus, it follows from the Theorem that only when $\alpha=0.5$ we have a nontrivial limiting point.

Let us briefly describe the structure of the paper.
In Section~\ref{s3} we introduce the necessary notation, describe the function spaces in which the results are obtained, and also give some information from probability theory and measure theory.
Sec\-tion~\ref{s4} briefly presents the results of the author's previous papers on the existence of a unique solution to the Cauchy problem  (see~\cite{1}) and the existence of a stationary measure (see~\cite{2}) for the system~\eqref{eq1.1}.
Section~\ref{s7} will be devoted to the formulation and
proof of the main result, which will be based on some integral properties of stationary measures for the system~\eqref{eq1.1}. This integral properties will be presented and proved in Section~\ref{s15}.
Section~\ref{s16} will contain a number of definitions and auxiliary results, on which the proof of the main result from \S\,\ref{s7} and the proof of integral properties from \S\,\ref{s15} will be based. Section~\ref{s11} will be devoted to the inviscid limit of stationary measures for one similar baroclinic atmosphere system.

\vskip15pt
\subsection*{Acknowledgments}

The author is grateful to Valentin~Dymnikov and Vladimir Krupchatnikov for posing
the problem, to Akif Ibragimov for  support, to Edriss S. Titi for inspiring discussion of the results of the work,
to Armen~Shirikyan and Vahagn~Nersesyan for helpful discussions.

\newpage
\section{Preliminaries}
\label{s3}

\subsection{Function spaces}
\label{s3.1}
Throughout, $X$ is a~Banach space with norm ${\|\cdot\|}_X$, $B_{X}(R)$ is the closed ball in~$X$ of radius~$R$, centred at the origin.

Let  $C(\mathbb R_{+};X)$ be the space of all continuous functions
$\psi \colon \mathbb R_{+} \to X$.
With the metric
$$
d(\psi,\theta)=\sum_{i=1}^\infty 2^{-i}
\frac{\max_{t\in [0,i]}\|\psi (t) - \theta (t)\|_X}{1+\max_{t\in [0,i]}\|\psi (t) - \theta (t)\|_X}
$$
the space $C(\mathbb R_{+};X)$ is a~complete metric space.

Let $C_0^\infty(S)$ be the space of infinitely smooth functions $\psi \colon S\to {\mathbb R}$ satisfying
$$
\int_S \psi\,dS=0.
$$
Consider the family of norms
$$
\|\psi\|_{p}=((-\Delta)^p \psi,\psi)^{1/2},
\qquad p\in {\mathbb Z},
$$
in the space $C_0^\infty(S)$,
where $(\,\cdot\,{,}\,\cdot\,)$~is the  inner product in~$L_2(S)$.
For each $p\in {\mathbb Z}$
let $h^p$ be the completion
of the space~$C_0^\infty(S)$ in the norm~${\|\cdot\|}_p$.
We note that
\begin{equation}
\label{eq3.130}
  \text{$h^p$ is densely and compactly embedded in~$h^q$ for $q<p$, where $p,q \in {\mathbb Z}$.}
\end{equation}

The operator~$-\Delta$ acting on the space $C_0^\infty (S)$ with inner product $(\,\cdot \,{,}\,\cdot\,)$ is positive definite. Hence, $-\Delta$ can be extended to a~self-adjoint operator \begin{equation}
\label{eq3.140}
-\Delta \colon D(-\Delta)=h^{p+2}\subset h^p\to h^p,
\qquad p\in {\mathbb Z},
\end{equation}
where $D(A)$ is the domain of the operator~$A$,
and
\begin{equation}
\label{eq3.150}
\|\psi\|_{p}=\bigl((-\Delta)^{p/2} \psi, (-\Delta)^{p/2}\psi\bigr)^{1/2}
\quad
\text{for any }\ \psi \in h^p,
\quad p\in {\mathbb Z}.
\end{equation}

The generalized Cauchy--Schwarz inequality holds
\begin{equation}
\label{eq3.1}
|(\psi,\theta)| \le \|\psi \|_{p}\,\|\theta \|_{-p}
\quad \text{for any }\ \psi \in h^{p},
\quad\theta \in h^{-p},
\quad p\in {\mathbb Z},
\end{equation}
where the inner product $(\,\cdot \,{,}\,\cdot\,)$ extends to $h^p\times h^{-p}$ as follows:
\begin{equation}
\label{eq3.2}
(\psi,\theta)=\bigl((-\Delta)^{p/2}\psi, (-\Delta)^{-p/2}\theta\bigr),
\qquad p\in {\mathbb Z}.
\end{equation}

Let
\begin{equation}
\label{eq3.1000}
|\mspace{-1mu}\|\psi |\mspace{-1mu}\|_{p}=\langle A_0^{p/2} \psi, A_0^{p/2}\psi\rangle ^{1/2},
\quad p\in {\mathbb Z},
\end{equation}
$$
A_0=\begin{pmatrix}
-\Delta & 0\\
0 & -\Delta 	
\end{pmatrix},
\quad
\langle \psi,\theta \rangle=\int_{S}\langle \psi,\theta \rangle_{{\mathbb R}^2} \, dS,
$$
denote the norm on the space $H^p:=h^p\times h^p$, $p\in {\mathbb Z}$,
where $\langle \psi,\theta \rangle_{{\mathbb R}^2}$~is the inner product on~${\mathbb R}^2$.

Note that the definition of the spaces $h^p$ and $H^p$ can be extended to all real $p$ via the completion of the space
$C_0^\infty(S)$ with preservation of the properties
\eqref{eq3.130},
\eqref{eq3.140}, \eqref{eq3.1}, \eqref{eq3.2} and the norms~\eqref{eq3.150},~\eqref{eq3.1000} in these spaces. In~\cite{30}, \S\,II,  it is described in detail
how to do this. To define these spaces in~\cite{30}, \S\,II,  the notion of $(-\Delta)^p$ is introduced, where $p\in {\mathbb R}$. In what follows, we will use such spaces too.

Let
$$
\varkappa := C({\mathbb R}_{+}; H^2)\cap L_{2, {\rm loc}} ({\mathbb R}_{+}; H^3).
$$

We also set
\begin{equation}
\label{eq3.3}
l^+_2:=\bigl\{\{b_i\}_{i=1}^\infty \in l_2\colon b_i \ge 0,\,i=1, 2,\dots\bigr\},
\qquad
\mathrm b :=\|\{b_i\}_{i=1}^\infty \|_{l_2}=\biggl(\sum_ {i=1}^{\infty} b_i^2\biggr)^{1/2}.
\end{equation}

\subsection{Necessary facts from probability theory and measure theory}
\label{s3.2}

Let $\Omega $ be an arbitrary set and let $F$ be a~$\sigma$-algebra of subsets of~$\Omega $.
A~measurable space will be denoted by $(\Omega, F)$. Let ${\mathbb P}$ be
a~probability measure on a~measurable space $(\Omega, F)$. In what follows,
it will be assumed that $(\Omega, F, {\mathbb P})$ is a~complete probability space.

We let  ${\rm N}_{(\Omega, F, {\mathbb P})}$ denote the family of  subsets
 ${\mathrm M}_1\subset \Omega$ for which there exists a~subset ${\mathrm M}_2\subset \Omega$,
${\mathrm M}_2\in F$, such that ${\mathrm M}_1\subset {\mathrm M}_2$ and
${\mathbb P}({\mathrm M}_2)=0$.

A~filtration $\{F_t\}_{t\ge 0}$ in $(\Omega, F, {\mathbb P})$ is called \textit{normal}
if $\{F_t\}_{t\ge 0}$ is a~cad (right-con\-tin\-uous) filtration such that
all ${\mathbb P}$-nullsets of the $\sigma $-algebra~$F$ lie in~$F_0$.

In what follows we shall assume that the filtration $\{F_t\}_{t\ge 0}$ is
normal.

Consider two arbitrary measurable spaces $(\Omega_1,\sigma_1)$ and $(\Omega_2,\sigma_2)$.
A~mapping $\psi \colon (\Omega_1,\sigma_1) \to (\Omega_2,\sigma_2)$ is called
\textit{$\sigma_1|\sigma_2$-measurable} if the full inverse image
$$
\psi ^{-1}({\mathrm M}):=\{\omega \in \Omega_1\colon \psi (\omega)\in {\mathrm M}\}
$$
of any ${\mathrm M}\in \sigma_2$
under~$\psi $ lies in~$\sigma_1$.
Sometimes,
when it is clear from the context what $(\Omega_2,\sigma_2)$ is,
we may simply say that $\psi$ is \textit{$\sigma_1$-measurable}.
 Let $\psi\colon(\Omega_1,\sigma_1)\to (\Omega_2,\sigma_2)$
 be a~$\sigma_1$-measurable mapping and let $G\subset\sigma_1$~be a family of subsets of $\Omega _1$.  We let $\sigma (G, \psi) \subset \sigma_1$ denote the smallest $\sigma$-algebra containing $G$ and all the
full inverse images of sets from~$\sigma_2$ under the $\psi$.

We denote the Borel $\sigma $-algebra of a~topological space~${\mathscr T}$ by $\sigma_{\mathscr T}$.

Let ${\mathscr P}(X)$ be the family of probability measures on a~measurable space $(X,\sigma_X)$.
We let $C_b(X)$ denote the space of bounded continuous functions $\psi\colon{X\!\to\!{\mathbb R}}$ and equip $C_b(X)$ with the norm
\begin{equation}
\label{eq20.10}
\|\psi \|_{\infty, X}=\sup_{\vartheta \in X} |\psi(\vartheta)|.
\end{equation}
By definition, a~measure $m \in {\mathscr P}(X)$ acts on a~function $\psi \in C_b(X)$ as follows:
$$
(\psi,m)=\int_X \psi\,dm.
$$
Now $(\psi,m)$, $\psi \in C_b(X)$, is a~continuous linear functional on the space $C_b(X)$. Consequently,
${\mathscr P}(X)$ can be looked upon as a~subset of the dual space~$C^*_b(X)$, and hence, it can be endowed with the weak-$*$ topology.

Let $\{u^\omega (t, \vartheta),\, t\ge 0\}_{\vartheta \in X}$ be a~family of
homogeneous Markov processes with respect to the filtration $\{F_t\}_{t\ge 0}$, satisfying the Feller property, and with the same
family of transition operators.
The values of these processes lie in the separable
Banach space~$X$ and obey the initial condition $u^\omega(0,\vartheta)=\vartheta $, $\vartheta \in X$.
We assume that the paths $u^\omega(t,\vartheta)$, $t\ge 0$,
are continuous for almost all $\omega \in \Omega$ for any $\vartheta \in X$ and
are continuous in $\vartheta \in X$
for almost all $\omega \in \Omega$ on a~finite interval of variation of~$t$.
Let $P(t,\vartheta,\Gamma)$ be the transition function,
$$
P(t,\vartheta,\Gamma)={\mathbb P} \{\omega \in \Omega\colon u^\omega(t,\vartheta)\in \Gamma \},
\qquad
\vartheta \in X,
\quad \Gamma \in \sigma_X,
\quad t\ge 0.
$$
Consider the semigroup $P^*_t$,
\begin{equation}
\label{eq3.11}
P^*_t m (\Gamma)=\int_{X} P(t,\vartheta,\Gamma) \,m (d \vartheta),
\qquad
 m \in {\mathscr P}(X),
 \quad \Gamma \in \sigma_{X},
 \quad t\ge 0.
\end{equation}
A~measure $m \in {\mathscr P}(X)$ is called a~\textit{stationary measure for the family} $\{u^\omega(t,\vartheta),\allowbreak t\ge0\}_{\vartheta \in X}$, if
\begin{equation}
\label{eq3.12}
P^*_tm=m,
\qquad t\ge 0.
\end{equation}

\section{The summary of results from~\cite{1} and~\cite{2}}
\label{s4}

Assume we are given $F_0$-measurable random vector function  $u^{\omega}_0(x)$, where $u_0\in H^2$ almost surely.
We augment the system~\eqref{eq1.1} with the initial condition
\begin{equation}
\label{eq4.1}
u(0)=u_0.
\end{equation}

Let $\{E_i\}_{i=1}^{\infty}$ be an orthonormal basis for the space $H^0$.
As a~random vector
function~$\eta $, we consider
\begin{equation}
\label{eq4.2}
\eta ^{\omega} (t, x)=\frac{\partial}{\partial t}\zeta ^{\omega}(t,x),
\qquad
\zeta ^{\omega}(t, x)=\sum_ {i=1}^{\infty} b_i \beta_i^{\omega} (t)E_i,
\end{equation}
where $\{b_i\}_{i=1}^\infty \in l^+_2$, and $\{\beta_i^{\omega}(t)\}_{i=1}^{\infty}$, $t\ge 0$,
is a~sequence of independent real
Brownian motions with respect to the filtration $\{F_t\}_{t\ge 0}$.
Hence, the series in~\eqref{eq4.2} lies in the space $C({\mathbb R}_+; H^0)$ for almost all $\omega \in \Omega $
(see the proof of Theorem 4.3 in~\cite{19}, Ch.~I, \S\,4.1, for example).

We give the following definition of a~solution  $u^{\omega}(t,x, u_0)$, $t\ge 0$, of the Cauchy problem
\eqref{eq1.1},~\eqref{eq4.1},~\eqref{eq4.2}.

\begin{definition}
\label{defin1}
An $H^2$-valued random process $u^{\omega}(t,x,u_0)$, $t\geqslant 0$, defined on the
probability space $(\Omega,F,{\mathbb P})$, is said to be a~\textit{solution of the Cauchy problem~\eqref{eq1.1}, \eqref{eq4.1}, \eqref{eq4.2}}
if it has the following properties:

a) the random process $u(t)$, $t \geqslant 0$, is adapted to the filtration $\{F_t\}_{t \geqslant 0}$ and almost all its paths lie in the space~$\varkappa$;

b) almost all the paths of the random process~$u(t)$, $t\geqslant 0$, satisfy the equality
$$
A_1u(t)+\int_{0}^t\bigl(\nu A_2u(s)+A_3u(s)+B(u(s))\bigr)\,ds=
A_1u_0+\zeta(t), \qquad
t\ge 0,
$$
which is considered in the space $C({\mathbb R}_+;H^{-1})$.
\end{definition}

According to Theorem~2 in~\cite{1}, \S\,3,
for any $\nu >0$, $\gamma,\rho, k_0, k_1 \ge 0$, $\{b_i\}_{i=1}^{\infty}\in l^+_2$ and any $F_0$-measurable random vector function
~$u_0$,
where $u_0\in H^2$ almost surely,
there exists a~unique solution to the Cauchy problem~\eqref{eq1.1},~\eqref{eq4.1},~\eqref{eq4.2},
which depends continuously on a~finite interval of variation of~$t$
on the initial data $u_0$ for almost all $\omega \in \Omega $.
 Theorem~1 in~\cite{2}, \S\,3,  tells us that
the solution to the Cauchy problem~\eqref{eq1.1},~\eqref{eq4.1},~\eqref{eq4.2} is a~homogeneous Markov process
with respect to the filtration $\{F_t\}_{t \ge 0}$ which has the Feller property. Moreover, for any initial data $u_0\in H^2$
this random process has the same family of transition operators. According to~\cite{2}, \S\,4 (see also~\cite{3}, \S\,4, formulae (4.4)--(4.6)
and Lemma 1),
if
\begin{equation}
\label{eq4.3}
k_0<{\mathfrak F}(\nu,\gamma,\rho, k_1),
\end{equation}
where ${\mathfrak F}>0$ is some real function of the arguments $\nu $, $\gamma $, $\rho $, $k_1$,
then
there exits a~stationary measure $m_\nu \in {\mathscr P}(H^{2})$ of the family $\{u^{\omega}(t,x,\vartheta), \,t\ge 0\}_{\vartheta \in H^2}$ of homogeneous Markov processes with the Feller property, that is determined by
the solutions of the Cauchy problem~\eqref{eq1.1},
$u(0)=\vartheta$, $\vartheta \in H^2$.
We note that the function ${\mathfrak F}(\nu,\gamma,\rho, k_1)$ is greater than
\begin{equation}
\label{eq21.21}
\min\left\{4k_1, \frac{4(2+\gamma)}{(2-\gamma)^2}(2k_1+\rho)\right\}.
\end{equation}
 Therefore, if $k'_0$ satisfies the inequality~\eqref{eq7.10} or $k_0$ does not exceed  the quantity~\eqref{eq21.21} (here the condition on the numerical parameters in~\eqref{eq21.30} may not be met), then $k_0$ also satisfies the inequality~\eqref{eq4.3}.

We retain for the family $\{u^{\omega}(t,x,\vartheta),\, t\ge 0\}_{\vartheta \in H^2}$ all the notation  in~\S\,\ref{s3.2}  introduced
for the family $\{u^{\omega}(t,\vartheta),\, t\ge 0\}_{\vartheta \in X}$.

\section{The inviscid limit of stationary measures for the system~\eqref{eq1.1}}
\label{s7}

\subsection{Main result}
\label{s7.1}

In this section we shall assume that
$$
b_i=\nu ^\alpha {b'}_i, \ \ i=1,2,\ldots,
$$
where the sequence  $\{{b'}_i\}_{i=1}^{\infty}$ is independent of $\nu$ and non-zero, $\alpha$ is an arbitrary real number.
We define
$$
\zeta '= \sum_ {i=1}^{\infty} {b'}_i \beta_i^{\omega} (t)E_i.
$$
Then
$$
\zeta = \nu ^\alpha \zeta '.
$$
Recall that $\delta_0$ is the Dirac measure.
The following result holds.
\begin{theorem}
\label{t1}
{\rm I.} For any  $\gamma,\rho, k_0, k_1\ge 0$ satisfying the inequality  $k_0$ does not exceed the quantity ~\eqref{eq21.21}, any $\{{b'}_i\}_{i=1}^{\infty}\,{\in}\, l^+_2$ the following statements hold.

(i) Let $\alpha = 0.5$.
Then there exists a limiting point $m$ for any sequence of stationary measures
$\{m_{\nu_n}\}_{n=1}^{\infty}$,
$\nu_n\underset{n \to +\infty}{\longrightarrow} 0$, for the system ~\eqref{eq1.1}
in the sense of weak-$*$ convergence in ${\mathscr P}(H^{3-\varepsilon})$, where $\varepsilon \in (0,1]$ is an arbitrary real number and $m$ satisfies the equality
\begin{equation}
\label{eq7.20}
m\left(H^3\right)=1.
\end{equation}
The upper bounds also hold for the limiting measure:
\begin{equation}
\label{eq17.23}
 \int_{H^2} \langle A_2\vartheta, A_1 \vartheta\rangle  \, m (d\vartheta )
\le  {\rm C}_1,
\end{equation}
\begin{equation}
\label{eq27.19}
\int_{H^2} \exp\left(\kappa  |\mspace{-1mu}\|A_1\vartheta|\mspace{-1mu}\|^2_{0}\right)\,m  (d\vartheta) \le
{\rm C}_2,
\end{equation}
where  $\kappa , C_1, C_2$ are some positive constants.

(ii) Let $\alpha > 0.5$.
Then the probability measure  $m_{\nu}$  converges to the Dirac measure $\delta_0$  in the sense of weak-$*$ convergence in ${\mathscr P}(H^3)$ as $\nu\to 0$.

{\rm II.} Let, in addition, the equality~\eqref{eq21.30} holds for $\rho, k_0, k_1$.
This means that the inequality ~\eqref{eq7.10} holds for $\gamma,\rho', k'_0, k'_1$.  Then
in the case (i) the limiting measure additionally satisfies the conservation law
\begin{equation}
\label{eq17.48}
\int_{H^2} |\mspace{-1mu}\|\vartheta|\mspace{-1mu}\|^2_{2} \, m (d\vartheta ) +
\int_{H^2} \langle  A'_3\vartheta,  \vartheta\rangle \, m (d\vartheta ) =  {\rm C}_3
\end{equation}
and the lower bound
\begin{equation}
\label{eq18.20}
 \int_{H^2} |\mspace{-1mu}\|\vartheta|\mspace{-1mu}\|^2_{1} \, m  (d\vartheta )
\ge C_4,
\end{equation}
where $ C_3, C_4$ are some positive constants,
$$
A'_3=\begin{pmatrix}
-k'_0\Delta & 2k'_0\Delta
\\
k'_0\Delta & -(2k'_0+k'_1+\gamma)\Delta+\rho ' 	I
\end{pmatrix}, \quad \langle  A'_3\vartheta,  \vartheta\rangle \ge 0 \quad \text{for any } \vartheta \in H^2.
$$
Here
$m \ne \delta _0$.

Moreover, one more statement also holds.

(iii) Let $\alpha < 0.5$.
Then the set of probability measures   $\{m_{\nu}\}_{\nu >0}$ has no accumula\-tion points  in the sense of weak-$*$ convergence in ${\mathscr P}(H^2)$ as $\nu\to 0$.

\end{theorem}

\begin{remark}
\label{r103}
The constants in the formulation of Theorem~\ref{t1}  have the form:
$$
{\rm C}_1= \frac{\mathrm {b'} ^2}{2},  \qquad C_2 = \exp\left(3+ \frac{{\mathrm b'}^2}{\left(\sup_{i\ge 1} {b'}_i\right)^{2}}\right),
\qquad \kappa =  \frac{2 }{2+\gamma}  \left(\sup _{i\ge 1} {b'}_i\right)^{-2},
$$
\begin{equation}
\label{eq7.75}
\mathrm {b'} :=\|\{{b'}_i\}_{i=1}^\infty \|_{l_2}=\biggl(\sum_ {i=1}^{\infty} {b'}_i^2\biggr)^{1/2},
\end{equation}
\begin{equation}
\label{eq39.1}
\quad {\rm C}_3=\frac{1}{2} \left(\sum_{i=1}^\infty {b'}_i^2   \langle A^{-1}_1 E_i, E_i \rangle\right),
\quad
C_4=\frac{{\rm C}^2_3}{{\rm C}_1+ 2 (4k'_0+k'_1+\gamma +\rho'){\rm C}_3}.
\end{equation}

\end{remark}
\begin{remark}
\label{r2}
An analogue of  Theorem ~\ref{t1} can be obtained in a more general setting, when the right-hand side of the system~\eqref{eq1.1}
   has the form
      $f+\eta$, where $f=\nu f'$, $f'\in H^{-1}$ doesn't depend on  $\nu$.
\end{remark}

The proof of Theorem~\ref{t1} will be given in the next subsection~\S\,\ref{s7.2}.
It is worth noting that the scheme of the proof is partial analogue to the arguments in~\cite{18}, \S\,5.2.1, and to the scheme of the proof of Theorem 5.2.17 in~\cite{18}, \S\,5.2.4,  which are given in~\cite{18} for the solution of the Navier-Stokes equation perturbed by the white noise with respect to variables~$x$
 on a~two-dimensional torus, provided that
both the solution and the right-hand side of the system are divergence-free vector
functions in the variables~$x$.

We can do the next step and study the measure ${\bf m_\nu}$, which is given by the distribution of the stationary solution of the system~\eqref{eq1.1} as an element of the space $\varkappa$  almost surely.
By a stationary solution we mean a
solution $u(t)$, $t\ge 0$, of the Cauchy problem~\eqref{eq1.1},~\eqref{eq4.1},~\eqref{eq4.2}, which is a stationary random process with distribution  equals to the probability measure $m _\nu$ for each fixed $t\ge 0$.
In order to study the inviscid limit of the measure ${\bf m_\nu}$ and obtain a connection between the limiting measure and the system~\eqref{eq1.1} in the case $\nu=0$, we need smoothness~\eqref{eq7.20} due to the nonlinear term in the system~\eqref{eq1.1}.
If we confine ourselves to studying the limit of the measure $m_\nu$, then we can additionally obtain a simpler result, but with less smoothness of the limiting measure. This result partly was motivated by question which was asked by Professor Titi during my talk on the seminar "Nonlinear PDEs". Namely, the following theorem takes place. We assume that measures from the space $H^2$ extend to the spaces $H^{2-\varepsilon}$, $\varepsilon \in (0,1]$, in the same way as described in \S\,\ref{s16.1}.

\begin{theorem}
\label{t111}
  For any $\gamma,\rho \ge 0$, $k_0, k_1>0$ satisfying the inequality~$k_0$ is less than the quantity~\eqref{eq21.21}, any $\{{b'}_i\}_{i=1}^{\infty}\,{\in}\, l^+_2$ the following statements hold.

(i) Let $\alpha = 0$.
Then there exists a limiting point $m \ne \delta_0$ for any sequence of stationary measures
$\{m_{\nu_n}\}_{n=1}^{\infty}$,
$\nu_n\underset{n \to +\infty}{\longrightarrow} 0$, for the system ~\eqref{eq1.1}
in the sense of weak-$*$ convergence in ${\mathscr P}(H^{2-\varepsilon})$, where $\varepsilon \in (0,1]$ is an arbitrary real number and $m$ satisfies the equality
$$
m\left(H^2\right)=1.
$$
The limiting measure has some integral properties:
$$
 \int_{H^2} |\mspace{-1mu}\|\vartheta|\mspace{-1mu}\|^2_{2}  \, m (d\vartheta )
\le  \widetilde{{\rm C}}_1,
$$
$$
\int_{H^2} \langle  A^{\nu =0}_3\vartheta,  \vartheta\rangle \, m (d\vartheta ) =  {\rm C}_3,
$$
$$
 \int_{H^2} |\mspace{-1mu}\|\vartheta|\mspace{-1mu}\|^2_{0} \, m  (d\vartheta )
\ge \widetilde{\rm C}_2,
$$
$$
\int_{H^2} \exp\left(\widetilde{\kappa}  |\mspace{-1mu}\|A_1\vartheta|\mspace{-1mu}\|^2_{0}\right)\,m  (d\vartheta) \le
\widetilde{\rm C}_3,
$$
 where $A_3^{\nu =0}$ is the operator $A_3$ for $\nu =0$, ${\rm C}_3$ is defined in~\eqref{eq39.1}, $\widetilde{\kappa}$,
$\widetilde{\rm C}_1$--$\widetilde{\rm C}_3$ are some positive constants.

(ii) Let $\alpha > 0$.
Then the probability measure  $m_{\nu}$  converges to the Dirac measure $\delta_0$  in the sense of weak-$*$ convergence in ${\mathscr P}(H^2)$ as $\nu\to 0$.

(iii) Let $\alpha < 0$.
Then the set of probability measures   $\{m_{\nu}\}_{\nu >0}$ has no accumulation points  in the sense of weak-$*$ convergence in ${\mathscr P}(H^1)$ as $\nu\to 0$.

\end{theorem}

Since $\gamma,\rho \ge 0$, $k_0, k_1>0$ satisfy the inequality $k_0$ is less than the quantity~\eqref{eq21.21}, then it is easy to show that there exits  $\alpha >0$ independent of  $\nu$ such that
$$
\langle A_3\psi,A_1\psi\rangle\ge \alpha\langle A_1\psi, A_0\psi\rangle
\quad\text{for any }\
\psi \in H^2 \quad\text{and } \ \nu \ge 0.
$$
Using this inequality, we can estimate the solution of the system~\eqref{eq1.1} not with the help of the operator $A_2$, as will be done in the proof of Theorem~\ref{t1}, but with the help of the operator $A_3$.
The rest of the proof of this Theorem is similar to the proof of Theorem~\ref{t1}.

\begin{remark}
\label{r209}
An analogue of  Theorem ~\ref{t111} can be obtained in a more general setting, when the right-hand side of the system~\eqref{eq1.1}
   has the form
      $f+\eta$, where  $f\in H^{0}$ doesn't depend on  $\nu$.
\end{remark}

\subsection{Proof of Theorem~\ref{t1}}
\label{s7.2}
(i)
Consider the stationary measure  $m_\nu \in {\mathscr P}(H^2)$ for the system~\eqref{eq1.1}.
We can find  the probability space $(\Omega , F, {\mathbb P})$
(independent of $\nu$) and the  set of $H^2$-valued random variables $\{\upsilon_\nu\}_{\nu >0}$ defined on $(\Omega , F, {\mathbb P})$ such that  the distribution of $\upsilon_\nu$ is equal  to  $m _\nu$ (see the proof of the existence such probability space and random variables in ~\cite{50}, Ch. I, \S\,4, Exercise 6).
By Chebyshev's inequality  for real random variable $\xi \geqslant 0$ (see the proof in~\cite{21}, Ch. II, \S\,6.6):
$$
{\mathbb P}\{\xi \geqslant c\}\leqslant \frac{E\xi}{c}
\quad \text{for any}\quad
c>0
$$
 we have the inequality
\begin{equation}
\label{eq7.40}
{\mathbb P} \left(|\mspace{-1mu}\|\upsilon_\nu|\mspace{-1mu}\|^2_{3}\ge c^2\right)
\le \frac{ E|\mspace{-1mu}\|\upsilon_\nu|\mspace{-1mu}\|^2_{3}}{c^2}
\quad \text{for any}\quad
c>0\,\,\text{and }\, \nu >0.
\end{equation}
Using the inequality~\eqref{eq7.23} from Theorem~\ref{l7.6} and~\eqref{eq21.33},
we get
\begin{equation}
\label{eq7.50}
 \int_{H^2} |\mspace{-1mu}\|\vartheta|\mspace{-1mu}\|^2_{3}  \, m _\nu (d\vartheta ) \le \int_{H^2} \langle A_2\vartheta, A_1 \vartheta\rangle  \, m _\nu (d\vartheta )
\le  \frac{{\rm b}^2}{2\nu}
\end{equation}
 for any $\nu >0$, where ${\rm b}$ is defined in~\eqref{eq3.3}.
 Since the expectation can be expressed in terms of the distribution
of a random variable
(see ~\cite{19}, Ch.~I, \S\,1.1, formula~(1.8), for example), we find from~\eqref{eq7.40} and~\eqref{eq7.50} that
$$
{\mathbb P} \left(|\mspace{-1mu}\|\upsilon_\nu|\mspace{-1mu}\|_{3}> c\right)
\le {\mathbb P} \left(|\mspace{-1mu}\|\upsilon_\nu |\mspace{-1mu}\|^2_{3}\ge c^2\right)
\le \frac{{\rm b}^2}{2c^2\nu}
$$
for any $c>0$ and  $\nu >0$. Recall that the distribution $\upsilon_\nu$ is equal to $m_\nu$. Then
 we obtain
\begin{equation}
\label{eq7.70}
m_\nu \left(B_{H^3}(c)\right)
\ge 1 - \frac{{\rm {b'}}^2}{2c^2}
\end{equation}
for any $c>0$ and  $\nu >0$,  where $\mathrm {b'}$ is defined in~\eqref{eq7.75}.
From the proof of Lemma 2.4.6  in~\cite{15}, Ch.~II, \S\,4, from Theorem 1 in~\cite{24}, Ch.~XVII, \S\,1.1,  and Theorem 3, III, in~\cite{24}, Ch.~XVI, \S\,3.1, it follows that $B_{H^3}(c)$~is a compact subset of~$H^{3-\varepsilon}$ for any $c>0$ and $\varepsilon>0$. So from ~\eqref{eq7.70} for fixed
 $\varepsilon>0$  for any $\delta>0$ there exits a compact set $K_{\delta, \varepsilon} \subset H^{3-\varepsilon}$ such that
$$
m_\nu \left(K_{\delta, \varepsilon}\right)\ge 1 - \delta \qquad \text{for any } \nu>0.
$$
Then, using Theorem 11.5.4, (I) and (II), from~\cite{70}, Ch.~XI, \S\,5,  we obtain that
any sequence  $\{m_{\nu_n}\}_{n=1}^{\infty}$,
$\nu_n\underset{n \to +\infty}{\longrightarrow} 0$, has the limiting point $m$ in the sense of weak-$*$ convergence in ${\mathscr P}(H^{3-\varepsilon})$
for any
 $\varepsilon\in (0,1]$.  Since $B_{H^3}(c)$ is a closed set in~$H^{3-\varepsilon}$ for any $c>0$ and $\varepsilon\in (0,1]$, from~\eqref{eq7.70}, Theorem 11.1.1, (a) and (c), in~\cite{70}, Ch.~XI, \S\,1,  we have
\begin{equation}
\label{eq7.80}
m \left(B_{H^3}(c)\right)
\ge 1 - \frac{{\rm b'}^2}{2c^2}
\end{equation}
for any $c>0$.
We note that $H^3$ can be represented as a countable union of embedding balls of the form  $B_{H^3}(c)$ for some $c>0$, which are closed sets, and hence, Borel sets in $H^{3-\varepsilon}$ for any $\varepsilon\in (0,1]$. Then, by definition of a $\sigma$-algebra, the space $H^3$~is a Borel set in the space $H^{3-\varepsilon}$ for any $\varepsilon\in (0,1]$. Now \eqref{eq7.20} readily follows from~\eqref{eq7.80} and Theorem, (2), in ~\cite{21}, Ch.~II, \S\,1.2.

We also note that  we can assume that $m $ is independent of $\varepsilon \in (0,1]$. We will now prove this. Let $\varepsilon \in (0,1)$.  Since any closed ball from $H^{3 - \varepsilon}$ is a closed set in  $H^{2}$ and ~\eqref{eq2.10} is the metric, using
the continuity of the embedding of $H^{3 - \varepsilon}$ in  $H^{2}$,
Theorem~\ref{t55},   the uniqueness of the limit,  Proposition 11.2.3 on the extension of bounded Lipschitz function from~\cite{70}, Ch.~XI, \S\,2, it follows that
$m ^{\varepsilon}=m^{1}$ on the Borel $\sigma$-algebra of the space $H^{3 - \varepsilon}$.

 Let us prove the integral properties. We will prove the estimate~\eqref{eq17.23}. By Theorem~\ref{l7.6}, from the estimate \eqref{eq7.23} we obtain
\begin{equation}
\label{eq22.10}
 \int_{H^2} \langle A_2\vartheta, A_1 \vartheta\rangle  \, m _{\nu_n}(d\vartheta )
\le  {\rm c}_1, \qquad n=1,2,\ldots ,
\end{equation}
 where ${\rm c}_1$ is defined in~\eqref{eq21.33}. Note that ${\rm c}_1={\rm C}_1$ if $\alpha = 0.5$. Let $N$~is an arbitrary natural number.
We  let $H_N$ denote the linear hull of the set of eigenvector functions  of the operator $A_0$ corresponding to  its $N$ eigenvalues  $\lambda _k= k(k+1)$, $k=1,2,\ldots, N$,  and let $P_N\colon H^{2} \to H_N$~denote the orthogonal projector onto the set $H_N$. We denote by $F_{r, N}\colon H^2 \to {\mathbb R}$  for $r >0$, $N\in {\mathbb N}$ the function:
$$
F_{r, N}(\vartheta) = \left\{
\begin{array}{ll}
\langle A_2P_N\vartheta, A_1 P_N\vartheta\rangle , & \langle A_2P_N\vartheta, A_1 P_N\vartheta\rangle \le r^2,\\
r^2,& \langle A_2P_N\vartheta, A_1 P_N\vartheta\rangle > r^2.
\end{array}
\right.
$$
Clearly,  as $r\to +\infty$ it follows that
$$
F_{r, N}(\vartheta) \uparrow \langle A_2P_N\vartheta, A_1 P_N\vartheta\rangle
$$
 for any $\vartheta \in H^3$. As $N\to +\infty$ it also follows that
$$
\langle A_2P_N\vartheta, A_1 P_N\vartheta\rangle \uparrow \langle A_2\vartheta, A_1 \vartheta\rangle
$$
 for any $\vartheta \in H^3$.
Therefore, for any $r>0$ and $N \in {\mathbb N}$ from~\eqref{eq22.10}
 we have the inequalities
$$
 \int_{H^2}   F_{r, N}(\vartheta) \, m _{\nu_n}(d\vartheta )
\le  {\rm C}_1, \qquad n=1,2,\ldots .
$$
Clearly, the function $F_{r, N}\in C_b(H^2)$. Then  the inequality
$$
 \int_{H^2}   F_{r, N}(\vartheta) \, m (d\vartheta )
\le  {\rm C}_1
$$
follows from the existence the limiting point  $m $ for the sequence
 $\{m _{\nu_n}\}_{n=1}^\infty$ in the sense of weak-$*$ convergence in ${\mathscr P}(H^2)$.
Using  Fatou's lemma   (see~\cite{21}, Ch.~II, \S\,6.3, Theorem 2,\,a), for example) to pass to the limit as $r \to +\infty$ in this inequality, we obtain
$$
 \int_{H^2}   \langle A_2P_N\vartheta, A_1 P_N\vartheta\rangle  \, m (d\vartheta )
\le  {\rm C}_1.
$$
Thus, using  Fatou's lemma to pass to the limit as $N \to +\infty$ and taking~\eqref{eq7.20} into account, this establish~\eqref{eq17.23}.
The proof of the estimate~\eqref{eq27.19} is analogous.  We should use~\eqref{eq37.19} and~\eqref{eq37.1} instead of~\eqref{eq7.23} and~\eqref{eq21.33}, respectively.
The proof is a little easier since the function
 $F_{r}\colon H^2 \to {\mathbb R}$  for $r >0$:
$$
F_{r}(\vartheta) = \left\{
\begin{array}{ll}
\exp\left(\kappa |\mspace{-1mu}\|A_1\vartheta|\mspace{-1mu}\|^2_{0}\right) , & |\mspace{-1mu}\|A_1\vartheta|\mspace{-1mu}\|_{0} \le r,\\
\exp\left(\kappa r^2\right),& |\mspace{-1mu}\|A_1\vartheta|\mspace{-1mu}\|_{0}> r,
\end{array}
\right.
$$
belongs to the space $C_b(H^2)$ and, therefore, in order to carry out the passage to the limit with respect to the viscous parameter, it is not necessary to add projectors to the function.

We will prove~\eqref{eq17.48}. By Theorem~\ref{l7.6}, from the equality \eqref{eq7.48} we obtain
\begin{equation}
\label{eq22.30}
\int_{H^2} |\mspace{-1mu}\|\vartheta|\mspace{-1mu}\|^2_{2} \, m _{\nu_n} (d\vartheta ) + \frac{1}{\nu }
\int_{H^2} \langle  A_3\vartheta,  \vartheta\rangle \, m _{\nu_n} (d\vartheta ) =  {\rm c}_2, \qquad n=1, 2, \ldots,
\end{equation}
 where ${\rm c}_2$ is defined in~\eqref{eq21.33}. Clearly,  ${\rm c}_2={\rm C}_3$ if $\alpha = 0.5$. By~\eqref{eq21.30},
$$
\frac{1}{\nu} A_3 = A'_3.
$$
From~\eqref{eq22.30} for any $R>0$ we have
$$
\int_{B_{H^2}(R)} |\mspace{-1mu}\|\vartheta|\mspace{-1mu}\|^2_{2} \, m _{\nu_n} (d\vartheta ) +
\int_{H^2 \setminus B_{H^2}(R)} |\mspace{-1mu}\|\vartheta|\mspace{-1mu}\|^2_{2} \, m _{\nu_n}  (d\vartheta)
$$
\begin{equation}
\label{eq32.14}
+
\int_{B_{H^2}(R)} \langle  A'_3\vartheta,  \vartheta\rangle \, m _{\nu_n} (d\vartheta ) +
\int_{H^2 \setminus B_{H^2}(R)} \langle  A'_3\vartheta,  \vartheta\rangle \, m _{\nu_n} (d\vartheta )
=  {\rm C}_3, \qquad n=1, 2, \ldots.
\end{equation}
Using the inequalities~\eqref{eq27.19} and
~\eqref{eq7.52},
the obvious inequality
$$
|\mspace{-1mu}\|\vartheta|\mspace{-1mu}\|_{2} \le |\mspace{-1mu}\|A_1\vartheta|\mspace{-1mu}\|_{0}\qquad \text{for any } \vartheta \in H^2
$$
 and Lemma~\ref{l16.36} (let $G(t) = \exp\left(\kappa t\right)$, $t\ge 0$,), we establishe that the family of random variables $\{|\mspace{-1mu}\|\upsilon_{\nu_n}|\mspace{-1mu}\|^2_{2}\}^\infty_{n=1}$ is uniformly integrable. Then for any $\delta >0$ there exists $R_\delta >0$ such that for any $R\ge R_\delta $
\begin{equation}
\label{eq42.1}
\int_{H^2 \setminus B_{H^2}(R)} |\mspace{-1mu}\|\vartheta|\mspace{-1mu}\|^2_{2} \, m _{\nu_n}  (d\vartheta) < \delta , \qquad n=1,2, \ldots.
\end{equation}
It is easy to show from the inequality
$$
\|\psi \|^2_p \le 2^{-(q-p)}\|\psi \|^2_q \quad \text{for any }\ \psi \in h^q,\quad p, q \in {\mathbb Z},\quad q>p
$$
(see the inequality (16) in~\cite{30}), that
\begin{equation}
\label{eq31.1}
\langle A'_3 \psi, \psi\rangle \le  K|\mspace{-1mu}\|\psi|\mspace{-1mu}\|^2_{1} \le \frac{K}{2}|\mspace{-1mu}\|\psi|\mspace{-1mu}\|^2_{2}
\end{equation}
for any $\psi \in H^2$, where
$$
K\colon = 4k'_0+k'_1+\gamma +\rho'.
$$
Then from~\eqref{eq42.1} we obtain
$$
\int_{H^2 \setminus B_{H^2}(R)} \langle  A'_3\vartheta,  \vartheta\rangle \, m _{\nu_n} (d\vartheta ) < \frac{2}{K}\delta .
$$
Now  it follows from ~\eqref{eq32.14} and ~\eqref{eq42.1} that for any $\delta >0$ there exits  $R_\delta >0$ such that for any  $R\ge R_\delta $ we have the chain of inequalities:
$$
{\rm C}_3 - \left(1+\frac{2}{K}\right)\delta \le \int_{B_{H^2}(R)} |\mspace{-1mu}\|\vartheta|\mspace{-1mu}\|^2_{2} \, m _{\nu_n} (d\vartheta )
+
\int_{B_{H^2}(R)} \langle  A'_3\vartheta,  \vartheta\rangle \, m _{\nu_n} (d\vartheta )
$$
\begin{equation}
\label{eq32.16}
\le \int_{H^2} |\mspace{-1mu}\|\vartheta|\mspace{-1mu}\|^2_{2} \, m _{\nu_n} (d\vartheta )
+
\int_{H^2} \langle  A'_3\vartheta,  \vartheta\rangle \, m _{\nu_n} (d\vartheta )
= {\rm C}_3, \qquad n=1, 2, \ldots.
\end{equation}
Note that the integrals  over the ball $B_{H^2}(R)$ in~\eqref{eq32.14} can be presented  as the integrals over the whole space $H^2$ with  integrands that vanishes outside the ball. The integrands is bounded but not continuous. By  Urysohn’s lemma,  there exits a continuous function such that  $0\le \theta_l \le 1$, $\theta_l = 1$ on the ball $B_{H^2}(R)$, $\theta _l = 0$ on the closure of the set  $H^2 \setminus B_{H^2}(R+l)$, where $l>0$~is an arbitrary real number.
Considering, instead of integrals over the ball from~\eqref{eq32.14}, integrals over the whole space $H^2$ with the same integrands multiplied by the function $\theta_l$, we obtain integrals with  bounded continuous integrands and can in them go to the limit with respect to the viscous parameter due to the existence of the limiting point $m $ of the sequence
  $\{m _{\nu_n}\}_{n=1}^\infty$ in the sense of weak-$*$ convergence in ${\mathscr P}(H^2)$.
Moreover, from the continuity of the probability measure and assertion (c) of Theorem 11.1.1 from~\cite{70}, Ch.~XI, \S\,1,
the difference between the original integrals and the corresponding integrals with $\theta_l$ tends to zero uniformly in the viscous parameter when $l\to 0$.
That's why
it is easy to show that for a subsequence $\{m_{\nu_{n_s}}\}^\infty_{s=1}$, weak-$*$ converging to $m$  in ${\mathscr P}(H^2)$ as $s \to \infty$, convergences for any $R>0$ take place:
$$
\int_{B_{H^2}(R)} |\mspace{-1mu}\|\vartheta|\mspace{-1mu}\|^2_{2} \, m _{\nu_{n_s}} (d\vartheta ) \to \int_{B_{H^2}(R)} |\mspace{-1mu}\|\vartheta|\mspace{-1mu}\|^2_{2} \, m  (d\vartheta ),
$$
$$
\int_{B_{H^2}(R)} \langle  A'_3\vartheta,  \vartheta\rangle \, m _{\nu_{n_s}} (d\vartheta ) \to \int_{B_{H^2}(R)} \langle  A'_3\vartheta,  \vartheta\rangle \, m  (d\vartheta )
$$
 as $s\to \infty $. Now  from~\eqref{eq32.16} we get~\eqref{eq17.48}, since $\delta $ is  arbitrary.

We will prove the estimate~\eqref{eq18.20}.
When the space $h^{3}$ is equipped with the inner product
$$
(\psi,\theta)_{3}=\bigl((-\Delta)^{3/2} \psi, (-\Delta)^{3/2}\theta \bigr)
$$
it becomes a Hilbert space.
Then, expanding $\psi$ in the orthonormal basis
$\{e_i\}^{\infty}_{i=1}$ for the Hilbert space $h^0$ ($h^3\subset h^0$) consisting of eigenfunctions of the Laplace--Beltrami operator,
it is easy to show that
$$
\|\psi \|^2_2 \le \|\psi \|_1 \|\psi \|_3 \quad \text{for any }\ \psi \in h^3.
$$
Therefore, by the Cauchy--Schwarz inequality, we have
\begin{equation}
\label{eq8.40}
 \int_{H^2} |\mspace{-1mu}\|\vartheta|\mspace{-1mu}\|^2_{2} \, m  (d\vartheta ) \le \sqrt{\int_{H^2} |\mspace{-1mu}\|\vartheta|\mspace{-1mu}\|^2_{1} \, m  (d\vartheta )}\sqrt{\int_{H^2} |\mspace{-1mu}\|\vartheta|\mspace{-1mu}\|^2_{3} \, m  (d\vartheta )}.
\end{equation}
Using~\eqref{eq17.48}, we find that
\begin{equation}
\label{eq8.50}
\int_{H^2} |\mspace{-1mu}\|\vartheta|\mspace{-1mu}\|^2_{2} \, m  (d\vartheta )
=  {\rm C}_3 - \int_{H^2} \langle  A'_3\vartheta,  \vartheta\rangle \, m (d\vartheta ).
\end{equation}
Then it follows from~\eqref{eq17.23},~\eqref{eq8.40},~\eqref{eq8.50} and~\eqref{eq31.1} that
$$
{\rm C}_3 -
K\int_{H^2} |\mspace{-1mu}\|\vartheta|\mspace{-1mu}\|^2_{1} \, m  (d\vartheta )
\le \sqrt{{\rm C}_1}\sqrt{\int_{H^2} |\mspace{-1mu}\|\vartheta|\mspace{-1mu}\|^2_{1} \, m  (d\vartheta )}.
$$
Next, examining the quadratic inequality with respect to
$\sqrt{\int_{H^2} |\mspace{-1mu}\|\vartheta|\mspace{-1mu}\|^2_{1} \, m  (d\vartheta )}$, it is easy to show~\eqref{eq18.20}. Note that from this inequality it follows that $m \ne \delta_0$, since the sequence  $\{{b'}_i\}_{i=1}^{\infty}$  is non-zero.

(ii) By the estimate~\eqref{eq7.23} of Theorem~\ref{l7.6} and~\eqref{eq21.33},
\begin{equation}
\label{eq17.10}
 \int_{H^2} \langle A_2\vartheta, A_1 \vartheta\rangle  \, m _\nu(d\vartheta )
\le  \frac{\mathrm {b} ^2}{2\nu} = \frac{\nu ^{2\alpha - 1}\mathrm {b'}^2}{2}.
\end{equation}
 Consider $\psi \in C_b(H^3)$ such that $\|\psi \|_{L, H^3}\le 1$, where $\|\cdot \|_{L, H^3}$ is defined in~\eqref{eq2.11}. Then, from~\eqref{eq17.40} we obtain
$$
\left|(\psi , m_\nu) - (\psi , \delta_0)\right|=\left|\int_{H^3}\psi (\vartheta) m _\nu(d\vartheta ) - \psi(0)\right| =
\left|\int_{H^3}(\psi (\vartheta) - \psi(0)) m _\nu(d\vartheta )\right|
$$
$$
 \le \int_{H^3}\left|\psi (\vartheta) - \psi(0)\right|m _\nu(d\vartheta ) \le  \int_{H^2} |\mspace{-1mu}\|\vartheta|\mspace{-1mu}\|_{3} \, m _\nu (d\vartheta ).
$$
Therefore, by the Cauchy--Schwarz inequality and~\eqref{eq17.10},
$$
\left|(\psi , m_\nu) - (\psi , \delta_0)\right| \le   \sqrt{\int_{H^2} |\mspace{-1mu}\|\vartheta|\mspace{-1mu}\|^2_{3} \, m _\nu (d\vartheta )}\le \frac{\nu ^{\alpha -0.5}\mathrm {b'}}{\sqrt{2}}.
$$
Using $\alpha> 0.5$ and \S\,\ref{s16.2}, this establishes that
 the probability measure $m_{\nu}$ converges to the Dirac measure  $\delta_0$  in the sense of weak-$*$ convergence in ${\mathscr P}(H^3)$ as $\nu\to 0$.

(iii)
In the proof of Theorem~\ref{l7.6}, after the formula~\eqref{eq7.24}, the existence of a probability space is proved and
solution $u(t)$, $t\ge 0$, of the Cauchy problem~\eqref{eq1.1}, \eqref{eq4.1}, \eqref{eq4.2} in this probability space such that it is a stationary random process with the distribution equal to the probability measure $m _\nu$ for each fixed $t\ge 0$.
We fix the constant  $\beta >0$ and consider the random process $v(t)\colon  = \nu ^\beta u (\nu ^\beta t)$, $t\ge 0$. It is easy to show that the random process $v(t)$, $t\ge 0$, satisfies the following equation:
$$
\frac{\partial}{\partial t}A_1 v+\nu^{1+\beta } A_2v+\nu^{1+\beta }A'_3v+B_\beta (v)=\nu ^{\alpha + \frac{3\beta }{2}}\frac{\partial}{\partial t}\widetilde{\zeta}(t,x),
\qquad t>0,
$$
where $\widetilde{\zeta}(t) \colon = \nu ^{-\frac{\beta}{2}} \zeta' (\nu ^\beta t)$, $t\ge 0$,
$$
B_\beta (v)=\bigl(J(\Delta v_1+2\nu^\beta \mu ,v_1)+J (\Delta v_2, v_2),\,
J(\Delta v_2 -\gamma v_2, v_1)+J(\Delta v_1+2\nu^\beta \mu , v_2)\bigr)^{\mathrm T}.
$$
Using the scaling properties of the Brownian motion, we can easy prove that the disribution of $\widetilde{\zeta}$ is the same as the distribution of $\zeta'$. Let $\beta\colon = \frac{1}{2} - \alpha $ and  $\widetilde{\nu}\colon  = \nu ^{1+\beta }$. Then the random process $v(t)$, $t\ge 0$, satisfies the following equation:
\begin{equation}
\label{eq17.76}
\frac{\partial}{\partial t}A_1 v+\widetilde{\nu} A_2v+\widetilde{\nu}A'_3v+B_{\beta} (v)=\sqrt{\widetilde{\nu}} \frac{\partial}{\partial t}\widetilde{\zeta}(t,x),
\qquad t>0.
\end{equation}

To prove that $\{m_\nu\}_{\nu >0}$ has no accumulation points as $\nu \to 0$, we argue by contradiction: suppose that there exits
the sequence $\{\nu_n\}^\infty_{n=1}$ which goes to zero such that  $\{m_{\nu_n}\}_{n=1}^\infty$
 weak-$*$ converges in  ${\mathscr P}(H^2)$.
By Theorem 11.5.4, (I) and (II),  in~\cite{70}, Ch.~XI, \S\, 5,  for any $\delta $ there exits the ball $B_{H^2}(R_\delta )$ such that for any  $R \ge R_\delta$
$$
m_{\nu_n} (H^2\setminus B_{H^2}(R ))={\mathbb P}\{ |\mspace{-1mu}\|u_{\nu_n}(0)|\mspace{-1mu}\|_{2} >R\}  \le \delta, \qquad n=1,2, \ldots.
$$
Therefore,
$$
{\mathbb P}\{ |\mspace{-1mu}\|u_{\nu_n}(0)|\mspace{-1mu}\|_{2} >R\}  \to 0 \qquad \text{ as } R \to \infty, \qquad  n=1,2, \ldots.
$$
 Then from the definition of the random process $v(t)$, $t\ge 0$, it is easy to show that for any $\delta >0$
\begin{equation}
\label{eq17.77}
{\mathbb P}\{ |\mspace{-1mu}\|v_{\nu_n}(0)|\mspace{-1mu}\|_{2} > \delta \} \to 0 \qquad \text{ as } \qquad n\to \infty.
\end{equation}
From  assertion (i) of  Theorem~\ref{t1} the sequence of stationary measures of the system~\eqref{eq17.76},
constructed by  $\widetilde{\nu}_n = \nu _n^{1+\beta}$, has the limiting point $\widetilde{m}$
in the sense of weak-$*$ convergence in ${\mathscr P}(H^2)$.
It follows from the inequality~\eqref{eq18.20}  that the limiting point $\widetilde{m}$ is not equal to the Dirac measure $\delta_0$.
Note that the equation~\eqref{eq17.76} differs from the equation~\eqref{eq1.1} by the nonlinear operator. However, it is easy to show that  the assertion (i) of Theorem~\ref{t1} also holds for such equation. From the existence of the limiting point  $\widetilde{m}$,~\eqref{eq17.77} and assertions (a), (b) of Theorem 11.1.1 in~\cite{70}, Ch.~XI, \S\,1,  we have that for any
$\delta >0$
$$
\widetilde{m}\left(H^2\setminus B_{H^2}(\delta)\right) = 0.
$$
Hence, $\widetilde{m} = \delta_0$. We arrive at a contradiction.

The proof of Theorem~\ref{t1} is complete.

\section{The integral properties of stationary measures for the system~\eqref{eq1.1}}
\label{s15}
In this section we will prove a number of inequalities and  one equality for some integrals with respect to stationary measures, which were used in the proof of Theorem~\ref{t1}.

\subsection{Formulation of the integral properties}
\label{s7.31}
In this section we shall assume that
the sequence $\{b_i\}_{i=1}^\infty $
is non-zero.
The following result holds.
\begin{theorem}
\label{l7.6}
For any  $\nu >0$, $\gamma,\rho, k_0, k_1\ge 0$ satisfying the inequality
$k_0$ does not exceed the quantity~\eqref{eq21.21} and
any $\{b_i\}_{i=1}^{\infty}\,{\in}\, l^+_2$
 a stationary measure  $m _\nu \in {\mathscr P}(H^2)$ for the system~\eqref{eq1.1} satisfies the inequalities
\begin{equation}
\label{eq7.23}
 \int_{H^2} \langle A_2\vartheta, A_1 \vartheta\rangle  \, m _\nu(d\vartheta )
\le  {\rm c}_1,
\end{equation}
\begin{equation}
\label{eq8.20}
 \int_{H^2} |\mspace{-1mu}\|\vartheta|\mspace{-1mu}\|^2_{1} \, m _\nu (d\vartheta )
\ge \frac{{\rm c}^2_2}{{\rm c}_1+ \frac{4}{\nu }{\rm c}_2 {\rm c}_3},
\end{equation}
\begin{equation}
\label{eq37.19}
\int_{H^2} \exp\left(\kappa _*\nu |\mspace{-1mu}\|A_1\vartheta|\mspace{-1mu}\|^2_{0}\right)\,m_\nu  (d\vartheta) \le
{\rm c}_4
\end{equation}
and the equality
\begin{equation}
\label{eq7.48}
\int_{H^2} |\mspace{-1mu}\|\vartheta|\mspace{-1mu}\|^2_{2} \, m _\nu (d\vartheta ) + \frac{1}{\nu }
\int_{H^2} \langle  A_3\vartheta,  \vartheta\rangle \, m _\nu (d\vartheta ) =  {\rm c}_2,
\end{equation}
where  $\kappa_*$, ${\rm c}_1$--${\rm c}_4$~are some positive constants.
In addition, the support of this statio\-nary measure for the system~\eqref{eq1.1} lies in the space $H^3$:
\begin{equation}
\label{eq17.40}
m_\nu\left(H^3\right)=1
\end{equation}
for each $\nu >0$.
\end{theorem}

\begin{remark}
\label{r310}
The constants in the formulation of  Theorem~\ref{l7.6}  have the form:
\begin{equation}
\label{eq21.33}
{\rm c}_1= \frac{\mathrm b ^2}{2\nu}, \qquad {\rm c}_2=\frac{1}{2\nu} \left(\sum_{i=1}^\infty b_i^2   \langle A^{-1}_1 E_i, E_i \rangle\right),
\qquad
{\rm c}_3=4k_0+k_1+\nu\gamma +\rho,
\end{equation}
\begin{equation}
\label{eq37.1}
c_4 = \exp\left(3+ \frac{{\mathrm b}^2}{\left(\sup_{i\ge 1} b_i\right)^{2}}\right),
\qquad \kappa_* =  \frac{2 }{2+\gamma}  \left(\sup _{i\ge 1} b_i\right)^{-2},
\end{equation}
where ${\mathrm b}$ is defined in~\eqref{eq3.3}.
\end{remark}

The proof of Theorem~\ref{l7.6} will be given a little later in~\S\,\ref{s7.100}.
In the next subsection we  formulate the auxiliary results, on which the proof of  Theorem~\ref{l7.6} will be based.

\subsection{Auxiliary results}
\label{s7.0}
  To prove Theorem~\ref{l7.6} we shall need the upper estimate from~\cite{52}, Lemma 6, for the quantity $E|\mspace{-1mu}\|A_1u(t)|\mspace{-1mu}\|_{0}$ for large $t$, where $u(t)$, $t\ge 0$, is the solution of
the Cauchy problem \eqref{eq1.1},  \eqref{eq4.1}, \eqref{eq4.2}.

\begin{lemma}
\label{l7.4}
Let $\alpha_1\in(0,2)$~be an arbitrary real number.
Then for any $\nu >0$, $\gamma,\rho, k_0, k_1 \ge 0$ satisfying
\begin{equation}
\label{eq7.2}
\langle A_3\psi,A_1\psi\rangle\ge 0
\quad\text{for any }\
\psi \in H^2,
\end{equation}
for any $\{b_i\}_{i=1}^{\infty}\,{\in}\, l^+_2$, 
any $F_0$-measurable random vector function
  $u_0 \in L_2(\Omega ; H^2)$
the solution $u(t)$, $t\ge 0$, of the Cauchy problem  \eqref{eq1.1},  \eqref{eq4.1}, \eqref{eq4.2} satisfies the inequality
\begin{align*}
\notag
&E |\mspace{-1mu}\|A_1u(t)|\mspace{-1mu}\|^2_{0}+((2 - \alpha_1)\nu)E\int_0^{t} e^{2\alpha_2 (s-t)}|\mspace{-1mu}\|u(s)|\mspace{-1mu}\|^2_{3} \, ds
\\
&\qquad
\le e^{-2\alpha_2 t} E|\mspace{-1mu}\|A_1 u_0|\mspace{-1mu}\|^2_{0}+\frac{1- e^{-2 \alpha_2 t}}{2\alpha_2}\mathrm b^2
\end{align*}
for $t\ge 0$, where ${\mathrm b}$ is defined in~\eqref{eq3.3} and $\alpha_2={4\alpha_1 \nu}/{(2+\gamma )^2}$.
\end{lemma}

 To prove Theorem~\ref{l7.6} we shall also need the upper estimate for
 the quantity
$$
E \exp(\kappa _* \nu |\mspace{-1mu}\|A_1u(t)|\mspace{-1mu}\|^2_{0})
$$
for large $t$, where $u(t)$, $t\ge 0$, is the solution of
the Cauchy problem \eqref{eq1.1},  \eqref{eq4.1}, \eqref{eq4.2} and $\kappa_*$ is defined in~\eqref{eq37.1}.

\begin{lemma}
\label{l37.4}
For any $\nu >0$, $\gamma,\rho, k_0, k_1\ge 0$ satisfying~\eqref{eq7.2},
for any $\{b_i\}_{i=1}^{\infty}\,{\in}\, l^+_2$,
any $F_0$-measurable random vector function 
  $u_0$, where $u_0\colon \Omega \to H^2$ is Bochner measurable function  and
$$
E\exp\left(\kappa_*\nu|\mspace{-1mu}\|A_1u_0|\mspace{-1mu}\|^2_{0}\right) <\infty
$$
(vector functions differing on ${\mathbb P}$-nullsets are regarded as identical),
 the solution $u(t)$, $t\ge 0$, of
the Cauchy problem \eqref{eq1.1}, \eqref{eq4.1}, \eqref{eq4.2} satisfies the inequality
$$
E \exp(\kappa _*\nu |\mspace{-1mu}\|A_1u(t)|\mspace{-1mu}\|^2_{0})
\le e^{-\alpha t} E\exp(\kappa _*\nu |\mspace{-1mu}\|A_1u_0|\mspace{-1mu}\|^2_{0})+\frac{C(\kappa_*, \alpha )}{\alpha}
$$
for $t\ge 0$,
where
$$
C(\kappa_*, \alpha) = \Bigr(\kappa _*\nu{\mathrm b}^2 + \alpha \Bigl)\exp\left(\frac{2+\gamma }{ 4\nu }\Bigr( \kappa _*\nu
{\mathrm b}^2 + \alpha\Bigl)\right),
$$
${\mathrm b}$ is defined in~\eqref{eq3.3}  and $\alpha>0$ is an arbitrary real number.
\end{lemma}

The proof of Lemma~\ref{l37.4} will be discussed a little later in~\S\,\ref{s7.803}.

We shall also need one equality for the solution $u(t)$, $t\ge 0$, of
the Cauchy problem \eqref{eq1.1},  \eqref{eq4.1}, \eqref{eq4.2}.

\begin{lemma}
\label{l7.8}
 For any $\nu >0$, $\gamma,\rho, k_0, k_1 \ge 0$ satisfying~\eqref{eq7.2}, any $F_0$-measurable random vector function
  $u_0 \in L_2(\Omega ; H^2)$
  the solution $u(t)$, $t\ge 0$, of
the Cauchy problem \eqref{eq1.1},  \eqref{eq4.1}, \eqref{eq4.2} satisfies the equality
$$
E\langle A_1 u(t),  u(t)\rangle + 2E \int _0^t\langle (\nu A_2 + A_3)u(s),  u(s)\rangle \, ds
$$
\begin{equation}
\label{eq7.25}
 = E\langle A_1 u_0,  u_0\rangle +   t\left(\sum_{i=1}^\infty b_i^2   \langle A^{-1}_1 E_i, E_i \rangle\right)
\end{equation}
for $t\ge 0$.
\end{lemma}

The proof of Lemma~\ref{l7.8} will be discussed a little later in~\S\,\ref{s7.800}.

\subsection{Proof of Theorem~\ref{l7.6}}
\label{s7.100}
We shall first show that the inequality
\begin{equation}
\label{eq7.52}
\int_{H^2} |\mspace{-1mu}\|A_1\vartheta|\mspace{-1mu}\|^2_{0} \, m_\nu (d\vartheta ) \le
\frac{(2+\gamma)^2\mathrm b^2}{16\nu} <\infty
\end{equation}
holds.
To do this we consider the function  $F_{r}\colon H^2 \to {\mathbb R}$ for $r >0$:
$$
F_{r}(\vartheta) = \left\{
\begin{array}{ll}
|\mspace{-1mu}\|A_1\vartheta|\mspace{-1mu}\|^2_{0}, & |\mspace{-1mu}\|A_1\vartheta|\mspace{-1mu}\|_{0}\le r,\\
r^2,& |\mspace{-1mu}\|A_1\vartheta|\mspace{-1mu}\|_{0}> r.
\end{array}
\right.
$$
From~\eqref{eq3.11} and~\eqref{eq3.12}, since $m_\nu$~is the stationary measure, we have
for any $t\ge0$ the equality
$$
\int_{H^2} F_{r}(\vartheta)\, m _\nu(d\vartheta)
= \int_{H^2}  F_{r}(\vartheta ) \left.\left(\int_{H^2} P(t, \upsilon,  \Gamma)\, m _\nu (d\upsilon)\right)\right|_{\Gamma = d\vartheta}.
$$
By the proof of Fubini's theorem in~\cite{36}, Ch. III, \S\,11.9,
$$
\int_{H^2} F_{r}(\vartheta)\, m _\nu (d\vartheta)=
\int_{H^2}  \left(\int_{H^2}  F_{r}(\vartheta ) P(t, \upsilon, d\vartheta )\right) \, m _\nu (d\upsilon).
$$
For convenience, we  replace $\vartheta$ by $\upsilon$ and $\upsilon$ by  $\vartheta$ in the right-hand side of the equality:
\begin{equation}
\label{eq7.21}
\int_{H^2} F_{r}(\vartheta)\, m _\nu (d\vartheta)=
\int_{H^2}  \left(\int_{H^2}  F_{r}(\upsilon) P(t, \vartheta, d\upsilon )\right) \, m _\nu (d\vartheta).
\end{equation}
We want to estimate the right-hand side of the equality.
By the definition of expectation
\begin{equation}
\label{eq7.1201}
\int_{H^2} F_{r}(\upsilon) P(t, \vartheta , d\upsilon ) \le E |\mspace{-1mu}\|A_1u^\omega (t,x,\vartheta)|\mspace{-1mu}\|^2_{0}.
\end{equation}
Note that the inequality~\eqref{eq4.3} with $\nu$  replaced by zero and the sign $<$  replaced by $\le $ has the form:
\begin{equation}
\label{eq301.1}
k_0 \le \inf_{i=1, 2, \dots} \frac{4\chi (j(i))}{(j(i)-\gamma)^2},
\end{equation}
$$
\chi(y)=k_1(y^2+\gamma y)+\rho(\gamma+y),\qquad
j(y)=y(y+1), \qquad
y \geqslant 0.
$$
It is easy to show that~\eqref{eq301.1} is equivalent to the inequality $k_0$ does not exceed the quantity~\eqref{eq21.21}.
It is also easy to show from the proof of Lemma 6 in~\cite{2}, \S\,4, that if the inequality~\eqref{eq301.1} holds then the inequality~\eqref{eq7.2} holds.
Then by  Lemma~\ref{l7.4} for any $\alpha_1\in (0,2)$  the inequality
\begin{align*}
\notag
&E |\mspace{-1mu}\|A_1u^\omega (t,x,\vartheta)|\mspace{-1mu}\|^2_{0}+((2 - \alpha_1)\nu)E\int_0^{t} e^{2\alpha_2 (s-t)}|\mspace{-1mu}\|u^\omega (s,x,\vartheta)|\mspace{-1mu}\|^2_{3} \, ds
\\
&\qquad
\le e^{-2\alpha_2 t} E|\mspace{-1mu}\|A_1 \vartheta|\mspace{-1mu}\|^2_{0}+\frac{1- e^{-2 \alpha_2 t}}{2\alpha_2}\mathrm b^2, \qquad t\ge 0,
\end{align*}
holds.
Now from~\eqref{eq7.1201} for
$|\mspace{-1mu}\|A_1\vartheta|\mspace{-1mu}\|_{0} \le R$ we obtain
$$
\int_{H^2} F_{r}(\upsilon) P(t, \vartheta, d\upsilon) \le
e^{-2\alpha_2 t} R^2+\frac{1- e^{-2 \alpha_2 t}}{2\alpha_2}\mathrm b^2, \qquad t\ge 0.
$$
Since $m _\nu$, $P(t, \vartheta, \cdot )$, $t\ge 0$,~are  probability measures, using the obvious inequality
$$
F_{r}(\vartheta) \le  r^2 \quad \text{for any }\  r>0,
$$
 we have
$$
\int_{H^2}  \left(\int_{H^2}  F_{r}(\upsilon) P(t, \vartheta, d\upsilon )\right) \, m _\nu (d\vartheta)
 \le
r^2m _\nu(\{\vartheta \in H^2\colon |\mspace{-1mu}\|A_1\vartheta|\mspace{-1mu}\|_{0}>R\})
$$
$$
+e^{-2\alpha_2 t} R^2+\frac{1- e^{-2 \alpha_2 t}}{2\alpha_2}\mathrm b^2, \qquad t\ge 0.
$$
We estimated the right-hand side of the equality~\eqref{eq7.21}.  From this estimate, the equality ~\eqref{eq7.21}, passing to the limit as $t\to+\infty$, we obtain the inequality
$$
\int_{H^2} F_{r}(\vartheta)\, m _\nu (d\vartheta)\le r^2m _\nu(\{\vartheta \in H^2\colon |\mspace{-1mu}\|A_1\vartheta|\mspace{-1mu}\|_{0}>R\})
+ \frac{1}{2\alpha_2}\mathrm b^2.
$$
Passing to the limit in this inequality as  $R\to +\infty$, we see that
$$
\int_{H^2} F_{r}(\vartheta)\, m _\nu (d\vartheta)\le
\frac{1}{2\alpha_2}\mathrm b^2.
$$
Clearly,  we can take the infimum with respect to $\alpha_2$, using that $\alpha_2={4\alpha_1 \nu}/{(2+\gamma )^2}$ and $\alpha_1\in (0,2)$~is an arbitrary real number. Then
$$
\int_{H^2} F_{r}(\vartheta)\, m _\nu (d\vartheta)\le
\frac{(2+\gamma)^2\mathrm b^2}{16\nu}.
$$
Using Fatou's lemma to pass to the limit as $r\to +\infty$,
we obtain~\eqref{eq7.52}.
The proof of the estimate~\eqref{eq37.19} is analogous. Here,
the function of $\vartheta$
$$
\left\{
\begin{array}{ll}
\exp\left(\kappa_* \nu |\mspace{-1mu}\|A_1\vartheta|\mspace{-1mu}\|^2_{0}\right), & |\mspace{-1mu}\|A_1\vartheta|\mspace{-1mu}\|_{0}\le r,\\
\exp\left(\kappa _*\nu r^2\right),& |\mspace{-1mu}\|A_1\vartheta|\mspace{-1mu}\|_{0}> r,
\end{array}
\right. \qquad r>0,
$$
should be taken for the function $F_{r}\colon H^2 \to {\mathbb R}$, $r >0$,
and  Lemma~\ref{l37.4} should be used instead of Lemma~\ref{l7.4}.

We will prove the estimate~\eqref{eq7.23}. As it was said above, the estimate ~\eqref{eq7.2} follows from the inequality
$k_0$ does not exceed the quantity~\eqref{eq21.21}. Then, by the Theorem  2 in~\cite{2}, \S\,4,  we get that the solution of the Cauchy problem~\eqref{eq1.1},~\eqref{eq4.1},~\eqref{eq4.2}   satisfies the inequality
$$
E|\mspace{-1mu}\|A_1u(t)|\mspace{-1mu}\|^2_{0}
 + 2\nu E \int _0^t \langle  A_2 u(s), A_1 u(s)\rangle  \, ds
$$
\begin{equation}
\label{eq7.24}
\le E|\mspace{-1mu}\|A_1u_0|\mspace{-1mu}\|^2_{0} + t{\rm b}^2
\end{equation}
 for $t\ge 0$, if $u_0 \in L_2(\Omega ; H^{2})$~is a $F_0$-measurable random vector function.
We can find  the probability space $(\Omega', F', {\mathbb P}')$
(independent of $\nu$) and  $H^2$-valued random variable $\upsilon_\nu$ defined on $(\Omega' , F', {\mathbb P}')$ such that  the distribution of $\upsilon_\nu$ is equal  to  $m_\nu$ (see the proof of the existence such probability space and random variable in ~\cite{50}, Ch. I, \S\,4, Exercise 6).
 Let us complete the probability space $(\Omega', F', {\mathbb P}')$ in the same way as described in \S\,\ref{s16.4}.
 We keep the same notation for the complete probability space.
Clearly,  the random variable $\upsilon_\nu$ has the same distribution $m _\nu$ on it.
  We take the direct product of the probability spaces     $(\Omega, F, {\mathbb P})$ and $(\Omega', F', {\mathbb P}')$ as the probability space  $(\mathbf{\Omega},\mathbf {F},\pmb {\mathbb P})$ (see the definition of the direct product of probability spaces in~\cite{49}, Ch.~I, \S\,4.1, for example).
We can assume that this direct product space is a~complete space.
If it is not, we complete~it in the same way as described in \S\,\ref{s16.4}.
We denote elements of $\mathbf{\Omega}$ by ${\boldsymbol \omega}=(\omega ,\omega')$, where $\omega \in \Omega$, $\omega' \in \Omega'$. Consider the Cauchy problem ~\eqref{eq1.1},~\eqref{eq4.1},~\eqref{eq4.2} on $(\mathbf{\Omega},\mathbf {F},\pmb {\mathbb P})$, where
$$
u_0^{\boldsymbol \omega}:= \upsilon^{\omega'},
\qquad \eta ^{\boldsymbol\omega} (t, x):=\frac{\partial}{\partial t}\zeta ^{\boldsymbol \omega}(t,x),
$$
$$
\zeta ^{\boldsymbol \omega}(t,x):= \sum_ {i=1}^{\infty} b_i \beta_i^{\boldsymbol \omega} (t)E_i := \sum_ {i=1}^{\infty} b_i \beta_i^{\omega} (t)E_i.
$$
Clearly, the random variable $u_0^{\boldsymbol \omega}$ is independent of the random variable  $\zeta ^{\boldsymbol \omega}$ by construction. Consider the following filtration on the space  $(\mathbf{\Omega},\mathbf {F},\pmb {\mathbb P})$:
$$
{\mathbf F}_t:={\mathrm N}_{(\mathbf{\Omega},\mathbf {F},\pmb {\mathbb P})}\cup \sigma \bigl(\{{\mathrm M}\times\Omega '\}_{{\mathrm M}\in F_t}, u_0^{\boldsymbol \omega}\bigr), \qquad t\ge 0
$$
(see the proof of the fact that ${\mathbf F}_t$ is a $\sigma$-algebra for any  $t\ge 0$ in ~\cite{25}, \S\,I.4, Proposition I.4.5, for example). Then all the hypotheses of Theorem 2 from~\cite{1}, \S\,3, on the existence of a unique solution of the Cauchy problem are satisfied.
By Theorem 1 from~\cite{2}, \S\,3, in particular, the solution of this Cauchy problem is a homogeneous Markov process. Therefore, from Proposition 1.4 in~\cite{28}, Ch.~III, \S\,1, the property 6 in~\cite{24}, \S\,4.2.2, and~\eqref{eq3.11},~\eqref{eq3.12} it follows that the solution of this Cauchy problem is a stationary random process  (see the definition in~\cite{28}, Ch.~I, \S\,3, Definition 3.4, for example) with  distribution   equals  the probability measure $m _\nu$ for each fixed $t\ge 0$.  Then from~\eqref{eq7.24}, taking~\eqref{eq7.52} into account,
by Fubini's theorem  (see \cite{23}, Ch.~VII, \S\,48, Remark, for example),
we obtain that for almost all $s\ge 0$ there exits $E\langle  A_2 u(s), A_1 u(s)\rangle$ and
the inequality
\begin{equation}
\label{eq81.30}
 \int _0^t E\langle  A_2 u(s), A_1 u(s)\rangle  \, ds
\le \frac{t}{2\nu} {\rm b}^2
\end{equation}
holds.
Note that $\sigma_{\mathbb R}$ coincides with the smallest  $\sigma$-algebra containing all the intervals of the form $(r_1,r_2]$, $r_1<r_2$ (see the proof in~\cite{21}, Ch.~II, \S\,2.2, for example). We will prove that the set
$$
{{\rm M}}_{r_1, r_2}=\{\vartheta \in H^3\colon r_1<|\mspace{-1mu}\|\vartheta|\mspace{-1mu}\|^2_3\le r_2\}
$$
for any $r_1$, $r_2$ satisfying $0\le r_1< r_2$ is a Borel subset of $H^2$.
 A complement to an open set is  closed  and vice versa, and hence, by the definition of a Borel  $\sigma$-algebra,  $\sigma _{H^2}$ coincides with the smallest  $\sigma$-algebra containing all closed subsets of $H^2$. Note that the equality
\begin{equation}
\label{eq8.30}
{\rm M}_{r_1, r_2}=\left((B_{H^3}(\sqrt{r_2}))^c\cup B_{H^3}(\sqrt{r}_1)\right)^c ,\qquad 0\le r_1< r_2,
\end{equation}
holds, where
$\cdot ^c={\mathrm M}\setminus \cdot $~is the complement to the subset $\cdot$ of the set ${\mathrm M}$, in our case ${\mathrm M}=H^2$.
From the proof of Lemma 2.4.6  in~\cite{15}, Ch.~II, \S\,4, Theorem 1 in~\cite{24}, \S\,17.1.1, and Theorem 3, III, in~\cite{24}, \S\,16.3.1,  it follows that $B_{{H}^3}(R)$~ is a compact subset of~${H}^2$ for any $R\ge 0$, and hence, closed and Borel. Therefore, by the definition of $\sigma$-algebra, from~\eqref{eq8.30} we obtain that  ${{\rm M}}_{r_1, r_2}$, $0\le r_1< r_2$,~is a Borel subset of $H^2$.
Then, by Theorem 1 in~\cite{24}, \S\,3.1.1,
 the function $|\mspace{-1mu}\|\cdot |\mspace{-1mu}\|^2_{3}\colon H^3 \subset H^2 \to {\mathbb R}$ is $\sigma_{H^2}$-measurable, and hence, the function
 $\langle  A_2 \cdot, A_1 \cdot\rangle \colon H^3 \subset H^2 \to {\mathbb R}$
is $\sigma_{H^2}$-measurable. In particular, $H^3$ is a Borel subset of $H^2$.
From the existence $E\langle  A_2 u(s), A_1 u(s)\rangle$ for almost all $s\ge 0$, ~\eqref{eq81.30} and Theorem 7 in~\cite{21}, Ch.~II, \S\,6.8,
it follows the inequality
$$
 \int_{H^3} \langle A_2\vartheta, A_1 \vartheta\rangle  \, m (d\vartheta )
\le  \frac{1}{2\nu} {\rm b}^2.
$$
By the definition of distribution, since almost all the paths of the random process $u(t)$, $t \ge 0$, lie in $\varkappa $,
it is easy to show~\eqref{eq17.40},
and hence,

$$
m_\nu(H^2 \setminus H^3)=0.
$$
Therefore, \eqref{eq7.23} holds.
A similar arguments show that from the equality~\eqref{eq7.25} we have~\eqref{eq7.48}.

The proof~\eqref{eq8.20} is analogous to the proof of the inequality~\eqref{eq18.20}, so we will not  give it.

The proof of Theorem~\ref{l7.6} is complete.

\subsection{Proof of Lemma~\ref{l7.8}}
\label{s7.800}
We will not  give the full proof of this Lemma, since it is analogous to the proof of Theorem 2 from~\cite{2}, \S\,4.
Note that in Lemma~5 from this  proof  we should take the spaces $H=H^{-1}$, $V=H^0$,
the functional  $\Phi(v)=\langle A_1^{-1}v,
 v\rangle $, $H^{-2}$  instead of $V^*$ and the It\^o process in $H^{-2}$ with constant diffusion  $A_1 u(t)$, $t\ge 0$ as the random process $\xi(t)$, $t\ge 0$.
Also, when constructing the proof, it is necessary to take  the equality
$$
\langle B(\psi), \psi \rangle = 0\quad\text{for any }\ \psi\in H^2
$$
into account, which takes place by Lemma 3 from~\cite{1}, \S\,3.
Inside the proof, the passage to the limit as $n \to \infty $ is carried out not by  Fatou's lemma, but by Lebesgue’s theorem on dominated
convergence
  (see~\cite{21}, Ch.~II, \S\,6.3, Theorem 3, for example).  This is because we are dealing with the proof of the equality, and not of the inequality, as in the proof of Theorem~2 from~\cite{2}, \S\,4.  To use Lebesgue’s theorem on dominated
convergence we should take account of the estimate
$$
E|\mspace{-1mu}\|A_1u(t)|\mspace{-1mu}\|^2_{0}
\le E|\mspace{-1mu}\|A_1u_0|\mspace{-1mu}\|^2_{0} + t{\rm b}^2
$$
for $t\ge 0$ (which is immediate corollary to~\eqref{eq7.24}).

\subsection{Proof of Lemma~\ref{l37.4}}
\label{s7.803}
We will not also give the full proof of this Lemma, since it is analogous to the proof of Lemma  6 in~\cite{52}. Here,  the functional \linebreak $e^{\alpha t}\exp (\kappa_* \nu |\mspace{-1mu}\|v|\mspace{-1mu}\|^2_{0})$ should be taken as the functional  $F (t, v)$.

\section{Appendix}
\label{s16}

\subsection{Extension of measure}
\label{s16.1}

 In this subsection we formulate the theorem on the extension of measure, which was proved in~\cite{113}, Ch. II,  \S\,2.1, Theorem 2.1 and its proof, and \S\,2.4, Remark 2.1.

We set
$$
\label{eq15.79}
\sigma_{\mathscr T}\cap {\mathrm M}:=\{{\mathrm N} \cap {\mathrm M}\colon {\mathrm N}\in \sigma_{\mathscr T}\}.
$$

\begin{theorem}
\label{t55}
 Let $\rm {X}_1$ and $\rm {X}_2$~be two  arbitrary metric spaces such that $\rm {X}_1$ is continuously embedded in $\rm {X}_2$.

{\rm I.} Then
\begin{equation}
\label{eq156.2}
\sigma _{\rm {X}_2}\cap\rm{X}_1 \subset \sigma _{\rm {X}_1}.
\end{equation}
Therefore, any measure ${m}\in {\mathscr P}(\rm {X}_1)$ can be extended to $\rm {X}_2$ as follows:
$$
m (\mathrm{M}) := m (\mathrm{M}\cap\mathrm{X}_1) \qquad \text{for any }\ \mathrm {M} \in  \sigma _{\mathrm {X}_2}.
$$

{\rm II.} In addition, if $\rm {X}_1$~is separable and any closed ball in $\rm {X}_1$ is closed in the topology of the space  $\rm {X}_2$, then
\begin{equation}
\label{eq156.3}
\sigma _{\rm {X}_1} \subset\sigma _{\rm {X}_2}, \qquad \text{in particular }\qquad \rm {X}_1 \in\sigma _{\rm {X}_2},
 \end{equation}
 and hence from~\eqref{eq156.2} and~\eqref{eq156.3}
$$
\sigma _{\rm {X}_2}\cap \rm {X}_1 = \sigma _{\rm {X}_1}.
$$
\end{theorem}

Thus,  when the measure is defined on $\sigma _{\rm {X}_1}$, one can always assume that it is defined on
$\sigma _{\rm {X}_2}$, if $\rm {X}_1$ is continuously embedded in  $\rm {X}_2$.

\subsection{Equivalence of topologies on the family of probability measures}
\label{s16.2}

Let Banach space $X$ be separable.
Then the weak-$*$ topology on~${\mathscr P}(X)$
is equivalent to the topology induced by the following distance between two measures:
\begin{equation}
\label{eq2.10}
\|m_1-m_2\|^*_{L, X}=\sup_{\substack{\psi \in C_b(X) \\
\|\psi \|_{L,X}\leqslant 1}}|(\psi,m_1)-(\psi,m_2)|, \qquad
m_1,m_2 \in {\mathscr P}(X)
\end{equation}
(the proof is given in~\cite{70}, Ch.~XI, \S\,3, Theorem 11.3.3, for example),
where
\begin{equation}
\label{eq2.11}
\|\psi\|_{L,X}=\|\psi\|_{\infty,X}+\sup_{\substack{\vartheta_1, \vartheta_2  \in X \\
\vartheta_1 \ne \vartheta_2}}
\frac{|\psi(\vartheta_1)-\psi(\vartheta_2)|}{\|\vartheta_1-\vartheta_2\|_X}\,,
\end{equation}
 the norm $\|\cdot\|_{\infty,X}$ is defined in~\eqref{eq20.10}.
In particular, the weak-$*$ convergence of a sequence of measures from~${\mathscr P}(X)$ to some limiting measure  is equivalent to the convergence to zero of the distance~\eqref{eq2.10} between the sequence terms and the limiting measure.

\subsection{Uniformly integrability of a family of random variables}
\label{s16.3}

In~\cite{21}, Ch.~II, \S\,6.4, Lemma~3, it was proved the following Lemma.

\begin{definition}
 A family   $\{\xi _n\}^\infty_{n=1}$ of random variables is said to be \textit {uniformly integrable}, if
$$
\sup_n E(|\xi_n| I_{|\xi_n| >R}) \to 0, \qquad R \to \infty.
$$
\end{definition}

\begin{lemma}
\label{l16.36}
 Let $\{\xi _n\}^\infty_{n=1}$~be a sequence of random variables with  finite expectation and $G(t)$, $t\ge 0$, be a nonnegative increasing function such that
$$
\lim_{t\to \infty} \frac{G(t)}{t}=\infty,
$$
$$
\sup_n EG(|\xi_n|) < \infty.
$$
Then the family  $\{\xi _n\}^\infty_{n=1}$ is uniformly integrable.
\end{lemma}

\subsection{The completion of a probability space}
\label{s16.4}
Let $(\Omega, F, {\mathbb P})$~ be some probabi\-lity space.
We complete $F$.
The \textit{completion}  of a~$\sigma $-algebra~$F$
with respect to the~probabi\-lity measure~${\mathbb P}$ is the smallest
~$\sigma $-algebra containing the sets of
$F \cup {\rm N}_{(\Omega, F, {\mathbb P})}$.
We point out  that the smallest ~$\sigma $-algebra containing some subsets of $\Omega$ always exits (for the proof see, for example,
\cite {21}, Ch.~II, \S\,2.1, Lemma 1).
Next, we extend the probability measure ${\mathbb P}$
on the completed $\sigma$-algebra by
${\mathbb P} ({\mathrm M}_1\cup {\mathrm M}_2):={\mathbb P} ({\mathrm M}_1)$, where ${\mathrm M}_1 \in  F$ and ${\mathrm M}_2 \in \mathrm N_{(\Omega, F, {\mathbb P})}$
(that this definition is consistent can be seen from~\cite{25}, \S\,I.4, Proposition I.4.5, for example).
Thus, we obtain the complete probability space.

\section{The inviscid limit of stationary measures for one similar  baroclinic atmosphere system}
\label{s11}

We note that a~slightly different
two-layer model for a~baroclinic atmosphere
on a~rotating two-dimensional unit sphere~$S$ with centre at the origin was considered
in~\cite{37}, Ch.~VI, \S\,1, and~\cite{38}.
When made dimensionless, the system of equations for this model differs
from~\eqref{eq1.1}
only in the operator $A_3$. Namely, in \cite{37}, Ch.~VI, \S\,1, and~\cite{38}
the operator~$A_3$ is as follows:
$$
\widehat{A}_3=\begin{pmatrix}
-k_0\Delta & k_0\Delta
\\
k_0\Delta & -(k_0+k_1)\Delta+\rho I
\end{pmatrix}.
$$
The existence of a~ stationary measure
for system~\eqref{eq1.1} with the operator~$A_3$ replaced by~$\widehat{A}_3$ was proved in~Theorem~4 from~\cite{2}, \S\,4.
Completely analogically to the arguments of the proofs of  Theorems~\ref{t1},~\ref{t111} and~\ref{l7.6}  we prove the following result.

\begin{theorem}
\label{t2}
The conclusions of Theorems~\ref{t1}, ~\ref{t111}, ~\ref{l7.6} and Remarks~\ref{r103}, \ref{r2}, \ref{r209}, \ref{r310} hold for the system ~\eqref{eq1.1} with  the operator~$A_3$ replaced by the operator~$\widehat{A}_3$,  taking into account  that
instead of
\begin{itemize}
\item the quantity ~\eqref{eq21.21} we have the quantity
$$
\frac{4(2+\gamma)}{\gamma^2}(2k_1+\rho),
$$
\item the right-hand side of the inequality~\eqref{eq7.10} we have the quantity
$$
\frac{4(2+\gamma)}{\gamma^2}(2k'_1+\rho'),
$$
\item the operator $A'_3$ we have the operator
$$
\widehat{A}'_3=\begin{pmatrix}
-k'_0\Delta & k'_0\Delta
\\
k'_0\Delta & -(k'_0+k'_1)\Delta+\rho ' 	I
\end{pmatrix}, \quad \langle  \widehat{A}'_3\vartheta,  \vartheta\rangle \ge 0 \quad \text{for any } \vartheta \in H^2,
$$
 \item the operator $A^{\nu =0 }_3$~we have the operator $\widehat{A}_3$,
 \item the constant ${\rm c}_3$~we have
the constant
${\rm \widehat{c}}_3 = 2k_0+k_1+\rho$,
\item the constant $C_4$~we have
 the constant
 $$
 \widehat{C}_4=\frac{{\rm C}^2_3}{{\rm C}_1+ 2 (2k'_0+k'_1+\rho'){\rm C}_3}.
 $$
 \end{itemize}
\end{theorem}

\end{fulltext}

\vskip10pt
\noindent
{\bf Yu. Yu. Klevtsova (Yulia Yu. Klevtsova)}
\\
Federal State Budgetary Institution ``Siberian Regional
\\
Hydrometeorological  Research  Institute'', Novosibirsk,
\\
Russian Federation;
\\
Siberian State University of Telecommunications and
\\
Information Sciences, Novosibirsk, Russian Federation
\\
{\it Email:} yy\_klevtsova@ngs.ru

\end{document}